\documentclass[a4paper,11pt]{amsart}

\usepackage{graphicx}
\usepackage{mathptmx}
\usepackage{amsmath}
\usepackage{amssymb}
\usepackage{enumitem}
\usepackage{xcolor} 
\usepackage{pgfplots}

\newmuskip\pFqmuskip

\newcommand*\pFq[6][8]{%
  \begingroup % only local assignments
  \pFqmuskip=#1mu\relax
  \mathcode`=\string"8000
  \begingroup\lccode`\~=`\,
  \lowercase{\endgroup\let~}\pFqcomma
  F^{#2}_{#3}{\left(\genfrac..{0pt}{}{#4}{#5}\bigg|#6\right)}%
  \endgroup
}
\newcommand{\pFqcomma}{\mskip\pFqmuskip}

\newtheorem{theorem}{Theorem}[section]
\newtheorem{lemma}[theorem]{Lemma}

\begin{document}

\title[]{Probabilistic heterogeneous Stirling numbers and Bell polynomials}

\author{Taekyun  Kim}
\address{Department of Mathematics, Kwangwoon University, Seoul 139-701, Republic of Korea}
\email{tkkim@kw.ac.kr}
\author{Dae San  Kim}
\address{Department of Mathematics, Sogang University, Seoul 121-742, Republic of Korea}
\email{dskim@sogang.ac.kr}


\begin{abstract}
Let $Y$ be a random variable satisfying specific moment conditions. This paper introduces and investigates probabilistic heterogeneous Stirling numbers of the second kind $H_{\lambda}^{Y}(n,k)$ and probabilistic heterogeneous Bell polynomials $H_{n,\lambda}^{Y}(x)$. These structures unify several classical and probabilistic families, including those of Stirling, Lah, Bell and Lah-Bell. By integrating the heterogeneous framework of Kim and Kim with probabilistic extensions, we derive explicit formulas, Dobi\'nski-like identities, and recurrence relations. We further establish connections to partial Bell polynomials and provide applications for Poisson and Bernoulli distributions.
\end{abstract}
 
\maketitle

\markboth{\centerline{\scriptsize Probabilistic heterogeneous Stirling numbers and Bell polynomials}}
{\centerline{\scriptsize Taekyun Kim and Dae San Kim}}

\section{Introduction} 
Let $Y$ be a random variable satisfying the moment conditions
\begin{equation}E\big[|Y|^{n}\big]<\infty,\quad n\in\mathbb{N}\cup{0},\quad \lim_{n\rightarrow\infty}\frac{|t|^{n}E[|Y|^{n}]}{n!}=0,\quad |t|<r, \label{1}\end{equation}
for some $r>0$, where $E$ denotes the mathematical expectation. Condition \eqref{1} ensures the existence of the moment generating function of $Y$, defined as
\begin{equation*}E\Big[e^{tY}\Big]=\sum_{n=0}^{\infty}\frac{t^{n}}{n!}E\big[Y^{n}\big], \quad (\text{see } [3,23,24]).
\end{equation*}
Furthermore, let $(Y_{j})_{j\ge 1}$ be a sequence of independent and identically distributed (i.i.d.) copies of $Y$, and let $S_{k}$ denote their partial sums:
\begin{displaymath}
S_{k}=\sum_{j=1}^k Y_j, \ (k\ge 1), \quad \text{with } S_{0}=0.
\end{displaymath} \par
Recently, Kim and Kim [18] introduced the heterogeneous Stirling numbers of the second kind, $H_{\lambda}(n,k)$, which generalize both the classical Stirling numbers of the second kind, ${n \brace k}$ (as $\lambda \rightarrow 0$), and the Lah numbers, $L(n,k)$ (as $\lambda \rightarrow 1$). They also investigated the heterogeneous Bell polynomials, $H_{n,\lambda}(x)$, which similarly generalize the Bell polynomials, $\phi_{n}(x)$, and the Lah-Bell polynomials, $LB_{n}(x)$. Special functions and polynomials serve as essential tools across mathematics, physics, and engineering. In recent years, significant research has focused on the degenerate versions and probabilistic extensions of these families. Degenerate polynomials are of interest not only for their number-theoretic and combinatorial properties but also for their applications in quantum mechanics and differential equations [6,14,17,29]. Concurrently, probabilistic extensions—including probabilistic Bell, Fubini, and Stirling polynomials—have emerged as a highly active area of study [13,15,16,19-21,24,26,27]. \par
This paper introduces and explores the probabilistic heterogeneous Stirling numbers of the second kind associated with $Y$, denoted by $H_{\lambda}^{Y}(n,k)$. These numbers generalize the probabilistic Stirling numbers of the second kind, ${n \brace k}_{Y}$, and the probabilistic Lah numbers, $L^{Y}(n,k)$. Additionally, we investigate the probabilistic heterogeneous Bell polynomials associated with $Y$, $H_{n,\lambda}^{Y}(x)$, which unify the probabilistic Bell polynomials, $\phi_{n}^{Y}(x)$, and the probabilistic Lah-Bell polynomials, $LB_{n}^{Y}(x)$. \par
The paper is organized as follows: Section 1 reviews the foundational definitions of Stirling, Lah, Bell, and Lah-Bell numbers and polynomials, alongside their probabilistic counterparts. It also recalls degenerate exponentials, heterogeneous Stirling numbers of the second kind associated with $Y$, heterogeneous Bell polynomials associated with $Y$ and the partial Bell polynomials. 
Section 2 presents our main results. We derive explicit expressions for $H_{\lambda}^{Y}(n,k)$ (Theorems 2.1, 2.2) and $L^{Y}(n,k)$ (Theorems 2.3, 2.4). For the polynomials $H_{n,\lambda}^{Y}(x)$, we provide explicit formulas (Theorems 2.5, 2.9, 2.11), a Dobiński-like formula (Theorem 2.6), recurrence relations (Theorems 2.7, 2.8), a binomial identity (Theorem 2.10), and derivative properties (Theorem 2.14). We also establish connections between $H_{\lambda}^{Y}(n,k)$ and partial Bell polynomials (Theorems 2.12, 2.13, 2.15). Finally, we provide applications for specific distributions: for a Poisson random variable, we relate the expectation of $\langle S_k \rangle_{n,\lambda}$ to a specific value of 
$H_{n,\lambda}(x)$ (Theorem 2.16) and get a sum representation for $H_{n,\lambda}^{Y}(x)$ involving the Bell polynomials (Theorem 2.17); for a Bernoulli random variable, we obtain explicit expressions for $H_{\lambda}^{Y}(n,k)$ and $H_{\lambda}^{Y}(n,k)$ (Theorem 2.18) and a sum representation involving the binomial coefficients for the expectation of $\langle S_{k}\rangle_{n,\lambda}$ (Theorem 2.20). General references for this study include [2,4,5,7,8,9,11,22,23,25]. \par
\vspace{0.1in}
The Stirling numbers of the second kind ${n \brace k},\ (n,k\ge 0)$, are defined by 
\begin{equation}
\frac{1}{k!}\big(e^{t}-1\big)^{k}=\sum_{n=k}^{\infty}{n \brace k}\frac{t^{n}}{n!},\quad (k\ge 0),\quad (\mathrm{see}\ [9,22]).\label{2}
\end{equation}
 The Lah numbers $L(n,k),\ (n,k\ge 0)$, are given by 
 \begin{equation}
 \frac{1}{k!}\bigg(\frac{1}{1-t}-1\bigg)^{k}=\sum_{n=k}^{\infty}L(n,k)\frac{t^{n}}{n!},\quad (k\ge 0). \label{3}
 \end{equation}
Thus, by \eqref{2} and \eqref{3}, we get 
\begin{equation}
{n \brace k}=\frac{1}{k!}\sum_{j=0}^{k}\binom{k}{j}(-1)^{k-j}j^{n},\quad (n\ge k), \label{4}
\end{equation}
and 
\begin{equation*}
L(n,k)=\frac{n!}{k!}\binom{n-1}{k-1},\quad (\mathrm{see}\ [9,10,12,22]). 
\end{equation*}
The probabilistic Stirling numbers of the second kind ${n \brace k}_{Y}$ and probabilistic Lah numbers $L^{Y}(n,k)$ are given by 
\begin{equation}
{n \brace k}_{Y}=\frac{1}{k!}\sum_{l=0}^{k}\binom{k}{l}(-1)^{k-l}E\big[S_{l}^{n}\big],\quad (\mathrm{see}\ [3,20]),\label{5}
\end{equation}
and 
\begin{equation}
L^{Y}(n,k)=\frac{1}{k!}\sum_{l=0}^{k}\binom{k}{l}(-1)^{k-l}E\Big[\langle S_{l}\rangle_{n}\Big],\quad (\mathrm{see}\ [10]),\label{6}
\end{equation}
where $\langle x\rangle_{0}=1,\ \langle x\rangle_{n}=x(x+1)\cdots(x+n-1),\ (n\ge 1)$. \par 
The Bell polynomials $\phi_{n}(x),\ (n\ge 0)$, are defined by 
\begin{equation}
e^{x(e^{t}-1)}=\sum_{n=0}^{\infty}\phi_{n}(x)\frac{t^{n}}{n!},\quad (\mathrm{see}\ [1,9,22]), \label{7}
\end{equation}
and Lah-Bell polynomials are given by 
\begin{equation}
e^{x\big(\frac{1}{1-t}-1\big)}=\sum_{n=0}^{\infty}\mathrm{LB}_{n}(x)\frac{t^{n}}{n!},\quad (\mathrm{see}\ [10,12]). \label{8}	
\end{equation}
From \eqref{7} and \eqref{8}, we have Dobinski-like formulas: 
\begin{equation}
\mathrm{LB}_{n}(x)=e^{-x}\sum_{k=0}^{\infty}\frac{\langle k\rangle_{n}}{k!}x^{k},\quad(n \ge 0),\quad (\mathrm{see}\ [10,12]),\label{9}	
\end{equation}
and 
\begin{equation}
\phi_{n}(x)=e^{-x}\sum_{k=0}^{\infty}\frac{k^{n}}{k!}x^{k},\quad (n\ge 0),\quad (\mathrm{see}\ [17,20-22,24]). \label{10}
\end{equation}
Recently, probabilistic Bell polynomials $\phi_{n}^{Y}(x)$, $(n\ge 0)$, and probabilistic Lah-Bell polynomials $\mathrm{LB}_{n}^{Y}(x),\ (n\ge 0)$, are given by 
\begin{equation}
\phi_{n}^{Y}(x)=e^{-x}\sum_{k=0}^{\infty}\frac{E[S_{k}^{n}]}{k!}x^{k},\quad (n\ge 0),\quad (\mathrm{see}\ [20,24]),\label{11}
\end{equation}
 and 
 \begin{equation}
 \mathrm{LB}_{n}^{Y}(x)=e^{-x}\sum_{k=0}^{n}\frac{E[\langle S_{k}\rangle_{n}]}{k!}x^{k},\quad (n\ge 0),\quad (\mathrm{see}\ [10,12]). \label{12}	
 \end{equation}
For any integer $k\ge 0$, the partial Bell polynomials are defined by 
\begin{equation}
\frac{1}{k!}\bigg(\sum_{m=1}^{\infty}x_{m}\frac{t^{m}}{m!}\bigg)^{k}=\sum_{n=k}^{\infty}B_{n,k}\big(x_{1},x_{2},\dots,x_{n-k+1}\big)\frac{t^{n}}{n!},\quad (\mathrm{see}\ [9,20,24]),\label{13}
\end{equation}
where 
\begin{align*}
&B_{n,k}(x_{1},x_{2},\cdots,x_{n-k+1})\\
&=\sum_{\substack{l_{1}+\cdots+l_{n-k+1}=k\\ l_{1}+2l_{2}+\cdots+(n-k+1)l_{n-k+1}=n } }\frac{n!}{l_{1}!l_{2}!\cdots l_{n-k+1}!}\bigg(\frac{x_{1}}{1!}\bigg)^{l_{1}} \bigg(\frac{x_{2}}{2!}\bigg)^{l_{2}}\cdots \bigg(\frac{x_{n-k+1}}{(n-k+1)!}\bigg)^{l_{n-k+1}},
\end{align*}
with $l_{1},l_{2},\dots,l_{n-k+1}$ nonnegative integers.
The complete Bell polynomials are given by 
\begin{equation}
\exp\bigg(\sum_{i=1}^{\infty}x_{i}\frac{t^{i}}{i!}\bigg)=\sum_{n=0}^{\infty}B_{n}\big(x_{1},x_{2},\cdots,x_{n}\big)\frac{t^{n}}{n!},\quad (\mathrm{see}\ [20,24]).\label{14}
\end{equation}
Thus, by \eqref{13} and \eqref{14}, we get 
\begin{equation*}
B_{n}(x_{1},x_{2},\cdots,x_{n})=\sum_{k=0}^{n}B_{n,k}(x_{1},x_{2},\cdots,x_{n-k+1}),\ (n\ge 0). 
\end{equation*} \par
For any nonzero $\lambda\in\mathbb{R}$, the degenerate exponentials are defined by 
\begin{equation}
e_{\lambda}^{x}(t)=\sum_{k=0}^{\infty}\frac{\langle x\rangle_{k,\lambda}}{k!}t^{k},\quad e_{\lambda}(t)=e_{\lambda}^{1}(t),\quad (\mathrm{see}\ [13-18]),\label{15}
\end{equation}
where 
\begin{equation*}
(x)_{0,\lambda}=1,\quad (x)_{n,\lambda}=x(x-\lambda)(x-2\lambda)\cdots\big(x-(n-1)\lambda\big),\ (n\ge 1). 
\end{equation*}
For any integer $n\ge 0$, the heterogeneous Stirling numbers of the second kind $H_{\lambda}(n,k),\ (n,k\ge 0)$, are defined by 
\begin{equation}
\langle x\rangle_{n,\lambda}=\sum_{k=0}^{n}H_{\lambda}(n,k)(x)_{k},\quad (\mathrm{see}\ [18]),\label{16}
\end{equation}
where $\langle x\rangle_{n,\lambda}=(-1)^{n}(-x)_{n,\lambda},\quad  (x)_{n}=(-1)^{n}\langle -x\rangle_{n}$. \\
From \eqref{15} and \eqref{16}, we note that 
\begin{equation}
\frac{1}{k!}\Big(e_{\lambda}^{-1}(-t)-1\Big)^{k}=\sum_{n=k}^{\infty}H_{\lambda}(n,k)\frac{t^{n}}{n!},\quad (\mathrm{see}\ [18]). \label{17}	
\end{equation}
Thus, by \eqref{17}, we get 
\begin{equation}
H_{\lambda}(n,k)=\frac{1}{k!}\sum_{j=0}^{k}\binom{k}{j}(-1)^{k-j}\langle j\rangle_{n,\lambda},\quad (n\ge k),\quad (\mathrm{see}\ [18]). \label{18}
\end{equation}
From \eqref{4}, \eqref{6} and \eqref{18}, we note that 
\begin{displaymath}
\lim_{\lambda\rightarrow 0}H_{\lambda}(n,k)={n \brace k}\quad\mathrm{and}\quad \lim_{\lambda\rightarrow 1}H_{\lambda}(n,k)=L(n,k). 
\end{displaymath}
The heterogeneous Bell polynomials, $H_{n,\lambda}(x),\ (n\ge 0)$, are defined by 
\begin{equation}
e^{x(e_{\lambda}^{-1}(-t)-1)}=\sum_{n=0}^{\infty}H_{n,\lambda}(x)\frac{t^{n}}{n!},\quad (\mathrm{see}\ [18]).\label{19}
\end{equation}
From \eqref{19}, we note that 
\begin{equation*}
H_{n,\lambda}(x)=\sum_{k=0}^{n}x^{k}H_{\lambda}(n,k),\ (n\ge 0),\quad (\mathrm{see}\ [18]), 
\end{equation*}
and 
\begin{equation*}
\lim_{\lambda\rightarrow 0}H_{n,\lambda}(x)=\phi_{n}(x),\quad\mathrm{and}\quad \lim_{\lambda\rightarrow 1}H_{n,\lambda}(x)=\mathrm{LB}_{n}(x). 
\end{equation*}
When $x=1,\ H_{n,\lambda}=H_{n,\lambda}(1)$ are called the heterogeneous numbers. \par 

\section{Probabilistic heterogeneous Stirling numbers of the second kind and Bell polynomials} 
Now, we consider a probabilistic extension of the heterogeneous Stirling numbers of the second kind, $H_{\lambda}^{Y}(n,k),\ (n,k\ge 0)$, which are called the {\it{probabilistic heterogeneous Stirling numbers of the second kind associated with $Y$}} and defined by 
\begin{equation}
\frac{1}{k!}\Big(E\Big[e_{\lambda}^{-Y}(-t)\Big]-1\Big)^{k}=\sum_{n=k}^{\infty}H_{\lambda}^{Y}(n,k)\frac{t^{n}}{n!},\quad (k\ge 0). \label{20}	
\end{equation}
From \eqref{20}, we have 
\begin{align}
&\sum_{n=k}^{\infty}H_{\lambda}^{Y}(n,k)\frac{t^{n}}{n!}=\frac{1}{k!}\Big(E\Big[e_{\lambda}^{-Y}(-t)\Big]-1\Big)^{k}\label{21}\\
&=\frac{1}{k!}\sum_{j=0}^{k}\binom{k}{j}(-1)^{k-j}\Big(E\Big[e_{\lambda}^{-Y}(-t)\Big]\Big)^{j}	=\frac{1}{k!}\sum_{j=0}^{k}\binom{k}{j}(-1)^{k-j}E\Big[e_{\lambda}^{-(Y_{1}+\cdots+Y_{j})}(-t)\Big] \nonumber\\
&=\frac{1}{k!}\sum_{j=0}^{k}\binom{k}{j}(-1)^{k-j}E\Big[e_{\lambda}^{-S}(-t)\Big]=\sum_{n=0}^{\infty}\frac{1}{k!}\sum_{j=0}^{k}\binom{k}{j}(-1)^{k-j}E\big[\langle S_{j}\rangle_{n,\lambda}\big]\frac{t^{n}}{n!}. \nonumber
\end{align}
Thus, by \eqref{21}, we get 
\begin{equation}
\frac{1}{k!}\sum_{j=0}^{k}\binom{k}{j}(-1)^{k-j}E\big[\langle S_{j}\rangle_{n,\lambda}\big]=\left\{\begin{array}{cc}
0, & \textrm{if $0 \le n<k$,}\\
H_{\lambda}^{Y}(n,k), & \textrm{if $n\ge k$.}  
\end{array}	\right.\label{22}
\end{equation}
Therefore, by \eqref{22}, we obtain the following theorem. 
\begin{theorem}
For $n\ge k\ge 0$, we have 
\begin{equation}
H_{\lambda}^{Y}(n,k)=\frac{1}{k!}\sum_{j=0}^{k}\binom{k}{j}(-1)^{k-j}E\big[\langle S_{j}\rangle_{n,\lambda}\big]. \label{23}
\end{equation}
\end{theorem}
In \eqref{23}, we have 
\begin{equation*}
\lim_{\lambda\rightarrow 1}H_{\lambda}^{Y}(n,k)=\frac{1}{k!}\sum_{j=0}^{k}\binom{k}{j}(-1)^{k-j}E\big[\langle S_{j}\rangle_{n}\big]=L^{Y}(n,k), 
\end{equation*}
and 
\begin{equation*}
\lim_{\lambda\rightarrow 0}H_{\lambda}^{Y}(n,k)=\frac{1}{k!}\sum_{j=0}^{k}\binom{k}{j}(-1)^{k-j}E\big[S_{j}^{n}\big]={n \brace k}_{Y}. 
\end{equation*}
From \eqref{20}, we note that 
\begin{align}
&\sum_{n=k}^{\infty}H_{\lambda}^{Y}(n,k)\frac{t^{n}}{n!}=\frac{1}{k!}\Big(E\Big[e_{\lambda}^{-Y}(-t)\Big]-1\Big)^{k}\label{24}\\
&=\frac{1}{k!}\Big(E\Big[e^{-\frac{Y}{\lambda}\log(1-\lambda t)}\Big]-1\Big)^{k}=\sum_{l=k}^{\infty}{l \brace k}_{Y}\frac{1}{l!}\bigg(-\frac{1}{\lambda}\log(1-\lambda t)\bigg)^{l}\nonumber\\
&=\sum_{l=k}^{\infty}{l \brace k}_{Y}\lambda^{-l}\frac{1}{l!}\bigg(\log\bigg(\frac{1}{1-\lambda t}\bigg)\bigg)^{l}=\sum_{l=k}^{\infty}{l \brace k}_{Y}\lambda^{-l}\sum_{n=l}^{\infty}{n \brack l}\lambda^{n}\frac{t^{n}}{n!} \nonumber\\
&=\sum_{n=k}^{\infty}\sum_{l=k}^{n}{l \brace k}_{Y}{n \brack l}\lambda^{n-l}\frac{t^{n}}{n!}, \nonumber
\end{align}
where ${n \brack k},\ (n,k\ge 0)$, are the unsigned Stirling numbers of the first kind defined by 
\begin{equation*}
\frac{1}{k!}\bigg(\log\bigg(\frac{1}{1-t}\bigg)\bigg)^{k}=\sum_{n=k}^{\infty}{n \brack k}\frac{t^{n}}{n!},\quad (k\ge 0),\quad (\mathrm{see}\ [18]). 
\end{equation*}
Therefore, by \eqref{24}, we obtain the following theorem. 
\begin{theorem}
For $n\ge k\ge 0$, we have 
\begin{equation*}
H_{\lambda}^{Y}(n,k)=\sum_{l=k}^{n}{l \brace k}_{Y}{n \brack l}\lambda^{n-l}. 
\end{equation*}
\end{theorem}
From \eqref{6}, we note that 
\begin{align}
&\sum_{n=k}^{\infty}L^{Y}(n,k)\frac{t^{n}}{n!}=\frac{1}{k!}\bigg(E\bigg[\bigg(\frac{1}{1-t}\bigg)^{Y}\bigg]-1\bigg)^{k}\label{25} \\
&=\frac{1}{k!}\Big(E\Big[e^{-Y\log(1-t)}\Big]-1\Big)^{k}=\sum_{l=k}^{\infty}{l \brace k}_{Y}\frac{1}{l!}\big(-\log(1-t)\big)^{l} \nonumber \\
&=\sum_{l=k}^{\infty}{l \brace k}_{Y}\sum_{n=l}^{\infty}{n \brack l}\frac{t^{n}}{n!}=\sum_{n=k}^{\infty}\sum_{l=k}^{n}{l \brace k}_{Y}{n \brack l}\frac{t^{n}}{n!}. \nonumber
\end{align}
Therefore, by \eqref{25}, we obtain the following theorem. 
\begin{theorem}
For $n\ge k\ge 0$, we have 
\begin{equation*}
L^{Y}(n,k)=\sum_{l=k}^{n}{l \brace k}_{Y}{n \brack l}. 
\end{equation*}
\end{theorem}
Let $\log_{\lambda}t$ be the compositional inverse function of $e_{\lambda}(t)$. Then we have 
\begin{equation*}
\log_{\lambda}(1+t)=\frac{1}{\lambda}\big((1+t)^{\lambda}-1\big),\quad (\mathrm{see}\ [14,16]). 
\end{equation*}
Note that 
\begin{equation*}
\log_{\lambda}\big(e_{\lambda}(t)\big)=e_{\lambda}(\log_{\lambda}(t))=t. 
\end{equation*}
From \eqref{20}, we have 
\begin{align}
&\sum_{n=k}^{\infty}L^{Y}(n,k)\frac{t^{n}}{n!}=\frac{1}{k!}\bigg(E\bigg[\bigg(\frac{1}{1-t}\bigg)^{Y}\bigg]-1\bigg)^{k}\label{26}\\
&=\frac{1}{k!}\bigg(E\Big[e_{\lambda}^{-Y}\big(\log_{\lambda}(1-t)\big)\Big]-1\bigg)^{k}=\sum_{l=k}^{\infty}H_{\lambda}^{Y}(l,k)\frac{1}{l!}\big(-\log_{\lambda}(1-t)\big)^{l} \nonumber\\
&=\sum_{l=k}^{\infty}H_{\lambda}^{Y}(l,k)\sum_{n=l}^{\infty}{n \brack l}_{\lambda}\frac{t^{n}}{n!}=\sum_{n=k}^{\infty}\sum_{l=k}^{n}H_{\lambda}^{Y}(n,k){n \brack l}_{\lambda}\frac{t^{n}}{n!}, \nonumber
\end{align}
where ${n \brack k}_{\lambda}$ are the degenerate unsigned Stirling numbers of the first kind defined by 
\begin{equation*}
\frac{1}{k!}\bigg(\log_{\lambda}\bigg(\frac{1}{1-t}\bigg)\bigg)^{k}=\frac{1}{k!}\Big(-\log_{\lambda}(1-t)\Big)^{k}=\sum_{n=k}^{\infty}{n \brack k}_{\lambda}\frac{t^{n}}{n!}, 
\end{equation*}
for $k\ge 0$, (see [14,15,16,18]). \\
Therefore, by \eqref{26}, we obtain the following theorem. 
\begin{theorem}
For $n\ge k\ge 0$, we have 
\begin{displaymath}
L^{Y}(n,k)=\sum_{l=k}^{n}H_{\lambda}^{Y}(l,k){n \brack l}_{\lambda}. 
\end{displaymath}
\end{theorem}
Now, we define the {\it{probabilistic heterogeneous Bell polynomials associated with $Y$}} by 
\begin{equation}
e^{x\big(E[e_{\lambda}^{-Y}(-t)] -1\big)}=\sum_{n=0}^{\infty}H_{n,\lambda}^{Y}(x)\frac{t^{n}}{n!}.\label{27}
\end{equation}
When $x=1,\ H_{n,\lambda}^{Y}=H_{n,\lambda}^{Y}(1),\ (n\ge 0)$, are called the {\it{probabilistic heterogeneous Bell numbers associated with $Y$}}. \\
By \eqref{27}, we get 
\begin{align}
&e^{x\big(E[e_{\lambda}^{-Y}(-t)] -1\big)}=\sum_{k=0}^{\infty}\frac{x^{k}}{k!}\Big(E\Big[e_{\lambda}^{-Y}(-t)\Big]-1\Big)^{k}\label{28} \\
&=\sum_{k=0}^{\infty}x^{k}\sum_{n=k}^{\infty}H_{\lambda}^{Y}(n,k)\frac{t^{n}}{n!}=\sum_{n=0}^{\infty}\sum_{k=0}^{n}x^{k}H_{\lambda}^{Y}(n,k)\frac{t^{n}}{n!}.\nonumber
\end{align}
Therefore, by \eqref{27} and \eqref{28}, we obtain the following theorem. 
\begin{theorem}
For $n\ge 0$, we have 
\begin{displaymath}
H_{n,\lambda}^{Y}(x)=\sum_{k=0}^{n}H_{\lambda}^{Y}(n,k)x^{k}. 
\end{displaymath}
\end{theorem}
From \eqref{27}, we note that 
\begin{align}
&\sum_{n=0}^{\infty}H_{n,\lambda}^{Y}(x)\frac{t^{n}}{n!}=e^{-x}e^{xE\big[e_{\lambda}^{-Y}(-t)\big] } \label{29} \\
&=e^{-x}\sum_{k=0}^{\infty}\frac{x^{k}}{k!}\Big(E\Big[e_{\lambda}^{-Y}(-t)\Big]\Big)^{k}=e^{-x}\sum_{k=0}^{\infty}\frac{x^{k}}{k!}E\Big[e_{\lambda}^{-S_{k}}(-t)\Big]\nonumber\\
&=e^{-x}\sum_{k=0}^{\infty}\frac{x^{k}}{k!}\sum_{n=0}^{\infty}E\Big[\langle S_{k}\rangle_{n,\lambda}\Big]\frac{t^{n}}{n!}=\sum_{n=0}^{\infty}e^{-x}\sum_{k=0}^{\infty}\frac{E[\langle S_{k}\rangle_{n,\lambda}]}{k!}x^{k}\frac{t^{n}}{n!}.\nonumber
\end{align}
Therefore, by \eqref{29}, we obtain the following Dobinski-like theorem. 
\begin{theorem}
For $n\ge 0$, we have 
\begin{equation*}
H_{n,\lambda}^{Y}(x)=e^{-x}\sum_{k=0}^{\infty}\frac{E\big[\langle S_{k}\rangle_{n,\lambda}\big]}{k!}x^{k}. 
\end{equation*}
\end{theorem} \par
Note that 
\begin{equation*}
\lim_{\lambda\rightarrow 1}H_{n,\lambda}^{Y}(x)=e^{-x}\sum_{k=0}^{\infty}\frac{E\big[\langle S_{k}
\rangle_{n}\big]}{k!}x^{k}=\mathrm{LB}_{n}^{Y}(x), 
\end{equation*}
and 
\begin{equation*}
\lim_{\lambda\rightarrow 0}H_{n,\lambda}^{Y}(x)=e^{-x}\sum_{k=0}^{\infty}\frac{E\big[S_{k}^{n}\big]}{k!}x^{k}=\phi_{n}^{Y}(x). 
\end{equation*}
By \eqref{27}, we get 
\begin{align}
\sum_{n=0}^{\infty}H_{n+1,\lambda}^{Y}(x)\frac{t^{n}}{n!}&=\frac{d}{dt}\sum_{n=0}^{\infty}H_{n,\lambda}^{Y}(x)\frac{t^{n}}{n!}=\frac{d}{dt} e^{x\big(E[e_{\lambda}^{-Y}(-t)] -1\big)}\label{30}\\
&=xE\Big[Ye_{\lambda}^{-(Y+\lambda)}(-t)\Big] e^{x\big(E[e_{\lambda}^{-Y}(-t)] -1\big)}\nonumber\\
&=x\sum_{k=0}^{\infty}E\big[\langle Y\rangle_{k+1,\lambda}\big]\frac{t^{k}}{k!}\sum_{m=0}^{\infty}H_{m,\lambda}^{Y}(x)\frac{t^{m}}{m!}\nonumber\\
&=\sum_{n=0}^{\infty}x\sum_{k=0}^{n}\binom{n}{k}E\big[\langle Y\rangle_{k+1,\lambda}\big]H_{n-k,\lambda}^{Y}(x)\frac{t^{n}}{n!}. \nonumber	
\end{align}
Therefore, by \eqref{30}, we obtain the following theorem. 
\begin{theorem}
For $n\ge 0$, we have 
\begin{equation*}
H_{n+1,\lambda}^{Y}(x)=x\sum_{k=0}^{n}\binom{n}{k}E\big[\langle Y\rangle_{k+1,\lambda}\big]H_{n-k,\lambda}^{Y}(x). 
\end{equation*}
\end{theorem}
Now, we observe that 
\begin{align}
&\sum_{m,n=0}^{\infty}H_{n+m,\lambda}^{Y}(t)\frac{x^{n}}{n!}\frac{y^{m}}{m!}=\sum_{m=0}^{\infty}\bigg(\frac{d}{dx}\bigg)^{m}\sum_{n=0}^{\infty}H_{n,\lambda}^{Y}(t)\frac{x^{n}}{n!}\frac{y^{m}}{m!} \label{31}\\
&=\sum_{m=0}^{\infty}\frac{y^{m}}{m!}\bigg(\frac{d}{dx}\bigg)^{m}\sum_{n=0}^{\infty}H_{n,\lambda}^{Y}(t)\frac{x^{n}}{n!}=e^{y\frac{d}{dx}} e^{t\big(E[e_{\lambda}^{-Y}(-x)] -1\big)}.\nonumber	
\end{align}
On the other hand, by binomial theorem, we get 
\begin{align}
&\sum_{m=0}^{\infty}\bigg(\frac{d}{dx}\bigg)^{m}\sum_{n=0}^{\infty}H_{n,\lambda}^{Y}(t)\frac{x^{n}}{n!}\frac{y^{m}}{m!}=\sum_{n=0}^{\infty}\frac{H_{n,\lambda}^{Y}(t)}{n!}\sum_{m=0}^{\infty}\bigg(\frac{d}{dx}\bigg)^{m}x^{n}\frac{y^{m}}{m!} \label{32} \\
&=\sum_{n=0}^{\infty}\frac{H_{n,\lambda}^{Y}(t)}{n!}(x+y)^{n}= e^{t\big(E[e_{\lambda}^{-Y}(-x-y)] -1\big)}.\nonumber
\end{align}
By \eqref{31} and \eqref{32}, we get 
\begin{align}
\sum_{m,n=0}^{\infty}H_{n+m,\lambda}^{Y}(t)\frac{x^{n}}{n!}\frac{y^{m}}{m!}&=e^{y\frac{d}{dx}} e^{t\big(E[e_{\lambda}^{-Y}(-x)] -1\big)}\label{33}\\
&= e^{t\big(E[e_{\lambda}^{-Y}(-x-y)] -1\big)}.\nonumber
\end{align}
From \eqref{15}, we note that 
\begin{align}
E\Big[e_{\lambda}^{-Y}(-x-y)\Big]&=E\bigg[(1-\lambda x)^{-\frac{Y}{\lambda}}\bigg(1-\frac{\lambda Y}{1-\lambda x}\bigg)^{-\frac{Y}{\lambda}}\bigg] \label{34}\\
&=E\bigg[e_{\lambda}^{-Y}(-x)e_{\lambda}^{-Y}\bigg(\frac{-y}{1-\lambda x}\bigg)\bigg]. \nonumber
\end{align}
By \eqref{33} and \eqref{34}, we get 
\begin{align}
&\sum_{m,n=0}^{\infty}H_{n+m,\lambda}^{Y}(t)\frac{x^{n}}{n!}\frac{y^{m}}{m!}=e^{t\big(E[e_{\lambda}^{-Y}(-x-y)] -1\big)}\label{35}\\
&= e^{t\big(E[e_{\lambda}^{-Y}(-x)] -1\big)}e^{t\big(E[e_{\lambda}^{-Y}(-x)]\big(e_{\lambda}^{-Y}\big(\frac{-y}{1-\lambda x}\big) -1\big)\big]\big)}\nonumber\\
&= e^{t\big(E[e_{\lambda}^{-Y}(-x)] -1\big)}\sum_{j=0}^{\infty}\frac{t^{j}}{j!}\bigg(E\bigg[e_{\lambda}^{-Y}(-x)\bigg(e_{\lambda}^{-Y}\bigg(\frac{-y}{1-\lambda x}\bigg)-1\bigg)\bigg]\bigg)^{j}. \nonumber
\end{align}
Note that 
\begin{align}
&\sum_{j=0}^{\infty}	\frac{t^{j}}{j!}\bigg(E\bigg[e_{\lambda}^{-Y}(-x)\bigg(e_{\lambda}^{-Y}\bigg(\frac{-y}{1-\lambda x}\bigg)-1\bigg)\bigg]\bigg)^{j}\label{36} \\
&=\sum_{j=0}^{\infty}\frac{t^{j}}{j!}E\bigg[e_{\lambda}^{-(Y_{1}+\cdots Y_{j})}(-x)\bigg(e_{\lambda}^{-Y_{1}}\bigg(\frac{-y}{1-\lambda x}\bigg)-1\bigg)\cdots\bigg(e_{\lambda}^{-Y_{j}}\bigg(\frac{-y}{1-\lambda x}\bigg)-1\bigg)\bigg] \nonumber \\
&=\sum_{j=0}^{\infty}\frac{t^{j}}{j!}\sum_{n=j}^{\infty}\bigg(\sum_{l_{1}+\cdots+l_{j}=n}\binom{n}{l_{1},\dots,l_{j}}E\Big[e_{\lambda}^{-S_{j}-n\lambda}(-x)\langle Y_{1}\rangle_{l_{1},\lambda}\cdots \langle Y_{j}\rangle_{l_{j},\lambda}\Big]\bigg)\frac{y^{n}}{n!}\nonumber\\
&=\sum_{n=0}^{\infty}\bigg(\sum_{j=0}^{n}\frac{t^{j}}{j!}\sum_{l_{1}+\cdots+l_{j}=n}\binom{n}{l_{1},\dots,l_{j}}E\bigg[e_{\lambda}^{-S_{j}-n\lambda}(-x)\prod_{i=1}^{j}\langle Y_{i}\rangle_{l_{i},\lambda}\bigg]\bigg)\frac{y^{n}}{n!}. \nonumber
\end{align}
From \eqref{35} and \eqref{36}, we have 
\begin{align}
&\sum_{m,n=0}^{\infty}H_{n+m,\lambda}^{Y}(t)\frac{x^{n}}{n!}\frac{y^{m}}{m!} 
= e^{t\big(E[e_{\lambda}^{-Y}(-x)] -1\big)}\label{37} \\
&\quad\quad\quad\quad\quad\times \sum_{j=0}^{\infty}\frac{t^{j}}{j!}\bigg(E\bigg[e_{\lambda}^{-Y}(-x)\bigg(e_{\lambda}^{-Y}\bigg(\frac{-y}{1-\lambda x}\bigg)-1\bigg)\bigg]\bigg)^{j} \nonumber\\
&=\sum_{n=0}^{\infty}\sum_{j=0}^{n}\frac{t^{j}}{j!}\sum_{l_{1}+\cdots+l_{j}=n}\binom{n}{l_{1},\dots,l_{j}}\frac{y^{n}}{n!} \nonumber\\
&\quad\quad\quad\quad\quad\times \sum_{l=0}^{\infty}E\bigg[\langle S_{j}+ny\rangle_{l,\lambda}\prod_{i=1}^{j}\langle Y_{i}\rangle_{l_{i},\lambda}\bigg]\frac{x^{l}}{l!}\sum_{k=0}^{\infty}H_{k,\lambda}^{Y}(t)\frac{x^{k}}{k!} \nonumber\\
&=\sum_{n=0}^{\infty}\sum_{m=0}^{\infty}\bigg(\sum_{j=0}^{n}\sum_{k=0}^{m}\binom{m}{k}\frac{t^{j}}{j!}\nonumber \\
&\quad\quad\quad\quad\quad\times \sum_{l_{1}+\cdots+l_{j}=n}\binom{n}{l_{1},\dots,l_{j}}E\bigg[\langle S_{j}+n\lambda\rangle_{m-k,\lambda}\prod_{i=1}^{j}\langle Y_{i}\rangle_{l_{i},\lambda}\bigg]\bigg)H_{k,\lambda}^{Y}(t)\frac{y^{n}}{n!}\frac{x^{m}}{m!}. \nonumber
\end{align}
Therefore, by comparing the coefficients on both sides of \eqref{37}, we obtain the following theorem. 
\begin{theorem}
For $m,n\ge 0$, we have 
\begin{align*}
H_{n+m,\lambda}^{Y}(t)&= \sum_{j=0}^{n}\sum_{k=0}^{m}\binom{m}{k}\frac{t^{j}}{j!}
\\
&\times\sum_{l_{1}+\cdots+l_{j}=n}\binom{n}{l_{1},\dots,l_{j}}E\bigg[\langle S_{j}+n\lambda\rangle_{m-k,\lambda}\prod_{i=1}^{j}\langle Y_{i}\rangle_{l_{i},\lambda}\bigg]H_{k,\lambda}^{Y}(t). 
\end{align*}
\end{theorem}
Taking $\lambda \rightarrow 0$ and $\lambda \rightarrow 1$ respectively gives
\begin{equation*}
\phi_{n+m}^{Y}(t)
=\sum_{j=0}^{n}\sum_{k=0}^{m}\binom{m}{k}\frac{t^{j}}{j!}\sum_{l_{1}+\cdots+l_{j}=n}\binom{n}{l_{1},\dots,l_{j}}E\bigg[S_{j}^{m-k}\prod_{i=1}^{j}Y_{i}^{l_{i}}\bigg]\phi_{k}^{Y}(t),
\end{equation*}
and 
\begin{equation*}
\mathrm{LB}_{n+m}^{Y}(t)=\sum_{j=0}^{n}\sum_{k=0}^{m}\binom{m}{k}\frac{t^{j}}{j!}\sum_{l_{1}+\cdots+l_{j}=n}\binom{n}{l_{1},\dots,l_{j}}E\bigg[\langle S_{j}+n\rangle_{m-k}\prod_{i=1}^{j}\langle Y_{i}\rangle_{l_{i}}\bigg]\mathrm{LB}_{k}^{Y}(t). 
\end{equation*}
From \eqref{27}, we have 
\begin{align}
&\sum_{n=0}^{\infty}H_{n,\lambda}^{Y}(x)\frac{t^{n}}{n!}=e^{x\big(E\big[e_{\lambda}^{-Y}(-t)\big]-1\big)} \label{38}\\
&=\exp\bigg(x\sum_{j=1}^{\infty}E\big[\langle Y\rangle_{j,\lambda}\big]\frac{t^{j}}{j!}\bigg)=\sum_{k=0}^{\infty}\frac{1}{k!}\bigg(x\sum_{j=1}^{\infty}E\big[\langle Y\rangle_{j,\lambda}\big]\frac{t^{j}}{j!}\bigg)^{k} \nonumber\\
&=\sum_{n=0}^{\infty}\sum_{k=0}^{n}x^{k} B_{n,k}\Big(xE\big[\langle Y\rangle_{1,\lambda}\big], xE\big[\langle Y\rangle_{2,\lambda}\big],\dots, xE\big[\langle Y\rangle_{n-k+1,\lambda}\big]\Big)\frac{t^{n}}{n!}. \nonumber
\end{align}
Therefore, by \eqref{38}, we obtain the following theorem. 
\begin{theorem}
For $n\ge 0$, we have 
\begin{equation*}
H_{n,\lambda}^{Y}(x)=\sum_{k=0}^{n}x^{k} B_{n,k}\Big(E\big[\langle Y\rangle_{1,\lambda}\big], E\big[\langle Y\rangle_{2,\lambda}\big],\dots, E\big[\langle Y\rangle_{n-k+1,\lambda}\big]\Big).
\end{equation*}
\end{theorem}
Note that 
\begin{align}
&\sum_{n=0}^{\infty}H_{n,\lambda}^{Y}(x+y)\frac{t^{n}}{n!}=e^{(x+y)(E[e_{\lambda}^{-Y}(-t)]-1)}  \label{39} \\
& e^{x(E[e_{\lambda}^{-Y}(-t)]-1)} e^{y(E[e_{\lambda}^{-Y}(-t)]-1)}=\sum_{k=0}^{\infty}H_{k,\lambda}^{Y}(x)\frac{t^{k}}{k!}\sum_{l=0}^{\infty}H_{l,\lambda}^{Y}(y)\frac{t^{l}}{l!} \nonumber\\
&=\sum_{n=0}^{\infty}\sum_{k=0}^{n}\binom{n}{k}H_{k,\lambda}^{Y}(x)H_{n-k,\lambda}^{Y}(y)\frac{t^{n}}{n!}.\nonumber
\end{align}
Therefore, by \eqref{39}, we obtain the following theorem. 
\begin{theorem}
For $n\ge 0$, we have 
\begin{displaymath}
H_{n,\lambda}^{Y}(x+y)=\sum_{k=0}^{n}\binom{n}{k}H_{k,\lambda}^{Y}(x)H_{n-k,\lambda}^{Y}(y). 
\end{displaymath}
\end{theorem}
By \eqref{13} and \eqref{27}, we get 
\begin{align}
&\sum_{n=0}^{\infty}H_{n,\lambda}^{Y}(x)\frac{t^{n}}{n!}= e^{x(E[e_{\lambda}^{-Y}(-t)]-1)}\label{40}\\
&=\Big(e^{E[e_{\lambda}^{-Y}(-t)]-1}-1+1\Big)^{x}=\sum_{k=0}^{\infty}\binom{x}{k}\Big(e^{E[e_{\lambda}^{-Y}(-t)]-1}-1\Big)^{k} \nonumber\\
&=\sum_{k=0}^{\infty}\binom{x}{k}k!\sum_{n=k}^{\infty}B_{n,k}\Big(H_{1,\lambda}^{Y},H_{2,\lambda}^{Y},\dots,H_{n-k+1,\lambda}^{Y}\Big)\frac{t^{n}}{n!} \nonumber\\
&=\sum_{n=0}^{\infty}\sum_{k=0}^{n}\binom{x}{k}k!B_{n,k}\Big(H_{1,\lambda}^{Y},H_{2,\lambda}^{Y},\dots,H_{n-k+1,\lambda}^{Y}\Big)\frac{t^{n}}{n!}. \nonumber
\end{align}
Therefore, by \eqref{40}, we obtain the following theorem. 
\begin{theorem}
For $n\ge 0$, we have 
\begin{equation*}
H_{n,\lambda}^{Y}(x)= \sum_{k=0}^{n}\binom{x}{k}k!B_{n,k}\Big(H_{1,\lambda}^{Y},H_{2,\lambda}^{Y},\dots,H_{n-k+1,\lambda}^{Y}\Big).
\end{equation*}
\end{theorem}
From \eqref{27}, we have 
\begin{equation}
te^{x(E[e_{\lambda}^{-Y}(-t)]-1)}=\sum_{j=1}^{\infty}jH_{j-1,\lambda}^{Y}(x)\frac{t^{j}}{j!}.\label{41}
\end{equation}
By \eqref{41}, we get 
\begin{align}
&\bigg(\sum_{j=1}^{\infty}jH_{j-1,\lambda}^{Y}(x)\frac{t^{j}}{j!}\bigg)^{k}= t^{k}e^{kx(E[e_{\lambda}^{-Y}(-t)]-1)} 
=t^{k}\sum_{j=0}^{\infty}k^{j}x^{j}\frac{1}{j!}\Big(E\big[e_{\lambda}^{-Y}(-t)\big]-1\Big)^{j} \label{42}\\
&=t^{k}\sum_{j=0}^{\infty}k^{j}x^{j}\sum_{n=j}^{\infty}H_{\lambda}^{Y}(n,j)\frac{t^{n}}{n!}=\sum_{n=0}^{\infty}\sum_{j=0}^{n}k^{j}x^{j}H_{\lambda}^{Y}(n,j)\frac{t^{n+k}}{n!} \nonumber\\
&=\sum_{n=k}^{\infty}k!\sum_{j=0}^{n-k}\binom{n}{k}k^{j}x^{j}H_{\lambda}^{Y}(n-k,j)\frac{t^{n}}{n!}.\nonumber
\end{align}
From \eqref{42}, we have 
\begin{align}
&\sum_{n=k}^{\infty}\sum_{j=0}^{n-k}\binom{n}{k}k^{j}x^{j}H_{\lambda}^{Y}(n-k,j)\frac{t^{n}}{n!}=\frac{1}{k!}\bigg(\sum_{j=1}^{\infty}jH_{j-1,\lambda}^{Y}(x)\frac{t^{j}}{j!}\bigg)^{k} \label{43} \\
&=\sum_{n=k}^{\infty}B_{n,k}\Big(H_{0,\lambda}^{Y}(x),2H_{1,\lambda}^{Y}(x),3H_{2,\lambda}^{Y}(x),\dots,(n-k+1)H_{n-k,\lambda}^{Y}(x)\Big)\frac{t^{n}}{n!}. \nonumber
\end{align}
Therefore, by \eqref{43}, we obtain the following theorem. 
\begin{theorem}
For $n\ge k\ge 0$, we have 
\begin{equation*}
\begin{aligned}
&\sum_{j=0}^{n-k}\binom{n}{k}k^{j}x^{j}H_{\lambda}^{Y}(n-k,j)\\
&= B_{n,k}\Big(H_{0,\lambda}^{Y}(x),2H_{1,\lambda}^{Y}(x),3H_{2,\lambda}^{Y}(x),\dots,(n-k+1)H_{n-k,\lambda}^{Y}(x)\Big).
\end{aligned}
\end{equation*}
\end{theorem}
Now, we observe that 
\begin{align}
&\sum_{n=k}^{\infty}B_{n,k}\Big(H_{1,\lambda}^{Y}(x),H_{2,\lambda}^{Y}(x),\dots,H_{n-k+1,\lambda}^{Y}(x)\Big)\frac{t^{n}}{n!} \label{44}	\\
&=\frac{1}{k!}\bigg(\sum_{j=1}^{\infty}H_{j,\lambda}^{Y}(x)\frac{t^{j}}{j!}\bigg)^{k}=\frac{1}{k!}\Big(e^{x(E[e_{\lambda}^{-Y}(-t)]-1)}-1\Big)^{k} \nonumber\\
&=\sum_{j=k}^{\infty}{j \brace k}x^{j}\Big(E\big[e_{\lambda}^{-Y}(-t)\big]-1\Big)^{j}=\sum_{j=k}^{\infty}{j \brace k}x^{j}\sum_{n=j}^{\infty}H_{\lambda}^{Y}(n,j)\frac{t^{n}}{n!} \nonumber\\
&=\sum_{n=k}^{\infty}\sum_{j=k}^{n}{j \brace k}H_{\lambda}^{Y}(n,j)x^{j}\frac{t^{n}}{n!}.\nonumber
\end{align}
Therefore, by \eqref{44}, we obtain the following theorem. 
\begin{theorem}
For $n\ge k\ge 0$, we have 
\begin{equation*}
B_{n,k}\Big(H_{1,\lambda}^{Y}(x),H_{2,\lambda}^{Y}(x),\dots,H_{n-k+1,\lambda}^{Y}(x)\Big)= \sum_{j=k}^{n}{j \brace k}H_{\lambda}^{Y}(n,j)x^{j}.
\end{equation*}
\end{theorem}
By \eqref{27}, we get 
\begin{align}
&\sum_{n=k}^{\infty}\bigg(\frac{d}{dx}\bigg)^{k}H_{n,\lambda}^{Y}(x)\frac{t^{n}}{n!}=\bigg(\frac{d}{dx}\bigg)^{k} e^{x(E[e_{\lambda}^{-Y}(-t)]-1)}\label{45}\\
&=\Big(E\big[e_{\lambda}^{-Y}(-t)\big]-1\Big)^{k} e^{x(E[e_{\lambda}^{-Y}(-t)]-1)}=k!\sum_{j=k}^{\infty}H_{\lambda}^{Y}(j,k)\frac{t^{j}}{j!}\sum_{i=0}^{\infty}H_{i,\lambda}^{Y}(x)\frac{t^{i}}{i!}\nonumber\\
&=\sum_{n=k}^{\infty}k!\sum_{i=0}^{n-k}\binom{n}{i}H_{i,\lambda}^{Y}(x)H_{\lambda}^{Y}(n-i,k)\frac{t^{n}}{n!}.\nonumber	
\end{align}
In particular, for $k=1$, then 
\begin{align}
\sum_{n=1}^{\infty}\frac{d}{dx}H_{n,\lambda}^{Y}(x)\frac{t^{n}}{n!}&=\Big(E\big[e_{\lambda}^{-Y}(-t)\big]-1\Big) e^{x(E[e_{\lambda}^{-Y}(-t)]-1)}\label{46} \\
&=\sum_{l=1}^{\infty}E\big[\langle Y\rangle_{l,\lambda}\big]\frac{t^{l}}{l!}\sum_{j=0}^{\infty}H_{j,\lambda}^{Y}(x)\frac{t^{j}}{j!} \nonumber\\
&=\sum_{n=1}^{\infty}\sum_{j=0}^{n-1}\binom{n}{j}E\big[\langle Y\rangle_{n-j,\lambda}\big]H_{j,\lambda}^{Y}(x)\frac{t^{n}}{n!}. \nonumber
\end{align}
Therefore, by \eqref{45} and \eqref{46}, we obtain the following theorem. 
\begin{theorem}
For $n\ge k\ge 1$, we have 
\begin{equation*}
\bigg(\frac{d}{dx}\bigg)^{k}H_{n,\lambda}^{Y}(x)=k!\sum_{j=0}^{n-k}\binom{n}{j}H_{j,\lambda}^{Y}(x)H_{\lambda}^{Y}(n-j,k), 
\end{equation*}
and 
\begin{equation*}
\bigg(\frac{d}{dx}\bigg)H_{n,\lambda}^{Y}(x)=\sum_{j=0}^{n-1}\binom{n}{j}E\big[\langle Y\rangle_{n-j,\lambda}\big]H_{j,\lambda}^{Y}(x). 
\end{equation*}
\end{theorem}
From \eqref{20}, we have 
\begin{align}
&\sum_{n=k}^{\infty}H_{\lambda}^{Y}(n,k)\frac{t^{n}}{n!}=\frac{1}{k!}\Big(E\big[e_{\lambda}^{-Y}(-t)\big]-1\Big)^{k}=\frac{1}{k!}\bigg(\sum_{j=1}^{\infty}E\big[\langle Y\rangle_{j,\lambda}\big]\frac{t^{j}}{j!}\bigg)^{k}\label{47}\\
&=\sum_{n=k}^{\infty}B_{n,k}\Big(E\big[\langle Y\rangle_{1,\lambda}\big], E\big[\langle Y\rangle_{2,\lambda}\big],\dots, E\big[\langle Y\rangle_{n-k+1,\lambda}\big]\Big)\frac{t^{n}}{n!}.\nonumber
\end{align}
Therefore, by \eqref{47}, we obtain the following theorem.
\begin{theorem}
For $n\ge k\ge 0$, we have 
\begin{displaymath}
H_{\lambda}^{Y}(n,k)= B_{n,k}\Big(E\big[\langle Y\rangle_{1,\lambda}\big], E\big[\langle Y\rangle_{2,\lambda}\big],\dots, E\big[\langle Y\rangle_{n-k+1,\lambda}\big]\Big).
\end{displaymath}
\end{theorem} 
Let $Y$ be the Poisson random variables with parameter $\alpha(>0)$ (see [23,28]). Then we have 
\begin{equation}
E\big[e_{\lambda}^{-Y}(-t)\big]=\sum_{n=0}^{\infty}e_{\lambda}^{-n}(-t)\frac{\alpha^{n}e^{-\alpha}}{n!}=e^{-\alpha}e^{\alpha e_{\lambda}^{-1}(-t)}=e^{\alpha(e_{\lambda}^{-1}(-t)-1)}. \label{48}	
\end{equation}
From \eqref{48}, we have 
\begin{align}
&\sum_{n=0}^{\infty}E\big[\langle S_{k}\rangle_{n,\lambda}\big]\frac{t^{n}}{n!}=E\big[e_{\lambda}^{-S_{k}}(-t)\big] \label{49} \\
&=E\Big[e_{\lambda}^{-(Y_{1}+\cdots+Y_{k})}(-t)\Big]=E\Big[e_{\lambda}^{-Y_{1}}(-t)\Big] E\Big[e_{\lambda}^{-Y_{2}}(-t)\Big]\cdots E\Big[e_{\lambda}^{-Y_{k}}(-t)\Big]\nonumber\\
&=e^{\alpha k(e_{\lambda}^{-1}(-t)-1)}=\sum_{n=0}^{\infty}H_{n,\lambda}(k\alpha)\frac{t^{n}}{n!}.\nonumber
\end{align}
Therefore, by \eqref{49}, we obtain the following theorem. 
\begin{theorem}
Let $Y$ be the Poisson random variable with parameter $\alpha(>0)$. Then we have 
\begin{align*}
E\big[\langle S_{k}\rangle_{n,\lambda}\big]=E\big[\langle Y_{1}+Y_{2}+\cdots+Y_{n}\rangle_{n,\lambda}\big]=H_{n,\lambda}(k\alpha),\ (n\ge 0). 
\end{align*}
\end{theorem}
Let $Y$ be the Poisson random variable with parameter $\alpha(>0)$. Then we have 
\begin{align}
&\sum_{n=0}^{\infty}H_{n,\lambda}^{Y}(x)\frac{t^{n}}{n!}=e^{x\big(E[e_{\lambda}^{-Y}(-t)]-1\big)}\label{50}\\
&=e^{x(e^{\alpha(e_{\lambda}^{-1}(-t)-1)}-1)}=\sum_{k=0}^{\infty}\phi_{k}(x)\frac{\alpha^{k}}{k!}\Big(e_{\lambda}^{-1}(-t)-1\Big)^{k} \nonumber\\
&=\sum_{k=0}^{\infty}\phi_{k}(x)\alpha^{k}\sum_{n=k}^{\infty}H_{\lambda}(n,k)\frac{t^{n}}{n!}=\sum_{n=0}^{\infty}\sum_{k=0}^{n}\phi_{k}(x)H_{\lambda}(n,k)\alpha^{k}\frac{t^{n}}{n!}. \nonumber
\end{align}
Therefore, by \eqref{50}, we obtain the following theorem. 
\begin{theorem}
Let $Y$ be the Poisson random variable with parameter $\alpha(>0)$. For $n\ge 0$, we have 
\begin{displaymath}
H_{n,\lambda}^{Y}(x)=\sum_{k=0}^{n}\phi_{k}(x)H_{\lambda}(n,k)\alpha^{k}. 
\end{displaymath}
\end{theorem}
Let $Y$ be the Bernoulli random variable with probability of success $p$ (see [23]). Then, for $n \ge 1$, we have 
\begin{equation}
E\big[\langle Y\rangle_{n,\lambda}\big]=\sum_{i=0}^{1}\langle i\rangle_{n,\lambda}p(i)=\langle 1\rangle_{n,\lambda}p.\label{51}
\end{equation}
From \eqref{20} and \eqref{51}, we note that 
\begin{align}
&\sum_{n=k}^{\infty}H_{\lambda}^{Y}(n,k)\frac{t^{n}}{n!}=\frac{1}{k!}\Big(E\big[e_{\lambda}^{-Y}(-t)\big]-1\Big)^{k}\label{52}\\
&=\frac{1}{k!}\bigg(\sum_{j=1}^{\infty}E\big[\langle Y\rangle_{j,\lambda}\big]\frac{t^{j}}{j!}\bigg)^{k}=\frac{p^{k}}{k!}\bigg(\sum_{j=1}^{\infty}\langle 1\rangle_{j,\lambda}\frac{t^{j}}{j!}\bigg)^{k} \nonumber\\
&=p^{k}\sum_{n=k}^{\infty}B_{n,k}\big(\langle 1\rangle_{1,\lambda},\langle 1\rangle_{2,\lambda},\dots,\langle 1\rangle_{n-k+1,\lambda}\big)\frac{t^{n}}{n!}.\nonumber	
\end{align}
Thus, by \eqref{52} and Theorem 2.15 with $Y=1$, we get 
\begin{equation}
H_{\lambda}^{Y}(n,k)=p^{k}B_{n,k}\big(\langle 1\rangle_{1,\lambda},\langle 1\rangle_{2,\lambda},\dots,\langle 1\rangle_{n-k+1,\lambda}\big)=p^{k}H_{\lambda}(n,k). \label{53}
\end{equation}
From \eqref{53} and Theorem 2.5 with $Y=1$, we have 
\begin{equation}
H_{n,\lambda}(xp)=\sum_{k=0}^{n}(px)^{k}H_{\lambda}(n,k)=\sum_{k=0}^{n}x^{k}H_{\lambda}^{Y}(n,k)=H_{n,\lambda}^{Y}(x). \label{54}	
\end{equation}
Therefore, by \eqref{54}, we obtain the following theorem. 
\begin{theorem}
Let $Y$ be the Bernoulli random variable with probability of success $p$. For $n\ge 0$, we have 
\begin{equation*}
H_{n,\lambda}^{Y}(x)=H_{n,\lambda}(xp),\quad H_{\lambda}^{Y}(n,k)=p^{k}H_{\lambda}(n,k). 
\end{equation*}
\end{theorem}
We observe that, for $n \ge 1$, we have
\begin{align}
&\langle 1\rangle_{n,\lambda}+\langle 2\rangle_{n,\lambda}+\cdots+\langle k\rangle_{n,\lambda}=\sum_{j=1}^{k}\langle j\rangle_{n,\lambda}\label{55}\\
&=\sum_{j=1}^{k}\sum_{l=1}^{n}H_{\lambda}(n,l)(j)_{l}=\sum_{l=1}^{n}H_{\lambda}(n,l)l!\sum_{j=1}^{k}\binom{j}{l}\nonumber\\
&=\sum_{l=1}^{n}H_{\lambda}(n,l)l!\sum_{j=1}^{k}\bigg(\binom{j+1}{l+1}-\binom{j}{l+1}\bigg)=\sum_{l=1}^{n}H_{\lambda}(n,l)l!\binom{k+1}{l+1}. \nonumber
\end{align}
Thus, by \eqref{55}, we obtain the following lemma. 
\begin{lemma}
For $k,n\in\mathbb{N}$, we have 
\begin{equation*}
\langle 1\rangle_{n,\lambda}+\langle 2\rangle_{n,\lambda}+\cdots+\langle k\rangle_{n,\lambda}=\sum_{l=1}^{n}H_{\lambda}(n,l)l!\binom{k+1}{l+1}. 
\end{equation*}
\end{lemma}
Let $Y$ be the Bernoulli random variable with probability of success $p$. Then we have 
\begin{align}
E\Big[e_{\lambda}^{-(Y_{1}+Y_{2}+\cdots+Y_{k})}\Big]&=E\big[e_{\lambda}^{-Y_{1}}(-t)\big]E\big[e_{\lambda}^{-Y_{2}}(-t)\big]\cdots E\big[e_{\lambda}^{-Y_{k}}(-t)\big] \label{56} \\
&=\Big(p\big(e_{\lambda}^{-1}(-t)-1\big)+1\Big)^{k},\nonumber
\end{align}
since
\begin{equation*}
E\big[e_{\lambda}^{-Y}(-t)\big]=(1-p)e_{\lambda}^{0}(-t)+pe_{\lambda}^{-1}(-t)=p\big(e_{\lambda}^{-1}(-t)-1\big)+1. 
\end{equation*}
By \eqref{56}, we get 
\begin{align}
&E\Big[e_{\lambda}^{-(Y_{1}+\cdots+Y_{k})}(-t)\Big]=\big(p(e_{\lambda}^{-1}(-t)-1)+1\big)^{k} \label{57} \\
&=\sum_{j=0}^{k}\binom{k}{j}p^{j}\big(e_{\lambda}^{-1}(-t)-1\big)^{j}=\sum_{j=0}^{k}\binom{k}{j}p^{j}j!\sum_{n=j}^{\infty}H_{\lambda}(n,j)\frac{t^{n}}{n!}\nonumber\\
&=\sum_{n=0}^{\infty}\sum_{j=0}^{n}\binom{k}{j}p^{j}j!H_{\lambda}(n,j)\frac{t^{n}}{n!}. \nonumber
\end{align}
On the other hand, by \eqref{17}, we get 
\begin{align}
E\Big[e_{\lambda}^{-(Y_{1}+Y_{2}+\cdots+Y_{k})}(-t)\Big]&=\sum_{n=0}^{\infty}E\big[\langle Y_{1}+\cdots+Y_{k}\rangle_{n,\lambda}\big]\frac{t^{n}}{n!} \label{58} \\
&=\sum_{n=0}^{\infty}E\big[\langle S_{k}\rangle_{n,\lambda}\big]\frac{t^{n}}{n!}. \nonumber
\end{align}

Therefore, by \eqref{57} and \eqref{58}, we obtain the following theorem. 
\begin{theorem}
Let $Y$ be the Bernoulli random variable with probability of success $p$. For $n\ge 0$, we have 
\begin{equation*}
E\big[\langle S_{k}\rangle_{n,\lambda}\big]=\sum_{j=1}^{n}\binom{k}{j}p^{j}j!H_{\lambda}(n,j). 
\end{equation*}
\end{theorem}
\section{Conclusion} 
In this study, we have successfully defined and analyzed the probabilistic heterogeneous Stirling numbers $H_{\lambda}^{Y}(n,k)$ and the probabilistic heterogeneous Bell polynomials $H_{n,\lambda}^{Y}(x)$. By leveraging the properties of the random variable $Y$ and its associated partial sums $S_k$, we established a diverse array of results, including: \\
$\indent \bullet$ Explicit Expressions: Theorems for both the numbers and polynomials. \\
$\indent \bullet$ Operational Identities: A Dobi\'nski-like formula and various recurrence relations. \\
$\indent \bullet$ Distributions: Specific applications involving Poisson and Bernoulli random variables. \\
These results demonstrate that the probabilistic heterogeneous framework provides a powerful tool for unifying diverse combinatorial identities. Future research may further explore the applications of these polynomials in the context of quantum mechanics or more complex stochastic processes.

\end{document}